\documentclass{amsart}
\usepackage{amssymb,amscd,amsmath,hyperref}
\usepackage{graphicx}
\theoremstyle{plain}
\newtheorem{theorem}[subsection]{Theorem}
\newtheorem{cor}[subsection]{Corollary}
\newtheorem{lem}[subsection]{Lemma}
\newtheorem{prop}[subsection]{Proposition}
\theoremstyle{definition}
\newtheorem{defn}[subsection]{Definition}
\theoremstyle{definition}
\newtheorem{rem}[subsection]{Remark}
\theoremstyle{definition}
\newtheorem{ex}[subsection]{Example}
\theoremstyle{definition}

\DeclareMathOperator{\N}{N}
\DeclareMathOperator{\Ne}{N}

\DeclareMathOperator{\ver}{ver}
\DeclareMathOperator{\Tor}{Tors}

\DeclareMathOperator{\supp}{supp}
\DeclareMathOperator{\Hom}{Hom}
\begin{document}
\title[The Laurent norm]%
{The Laurent norm}

\author[David G.~ Long]{David G. Long}
\address{Department of Mathematics, Northeastern University,
Boston, MA, 02115}
\email{\href{mailto:dlong2002us@yahoo.com}{dlong2002us@yahoo.com}}



\begin{abstract}
We generalize a semi-norm for the Alexander polynomial of a connected, compact, oriented 3-manifold on its first cohomology group to a semi-norm for an arbitrary Laurent polynomial $f$ on the dual vector space to the space of exponents of $f$.  We determine a decomposition formula for this Laurent norm; an expression for the Laurent norm for $f$ in terms of the Laurent norms for each of the irreducible factors of $f$.  For an $n$-variable polynomial $f$, we introduce a space of \( m \leq n\) essential variables which determine the reduced Laurent norm unit ball; a convex polyhedron of the same dimension $m$ as the Newton polyhedron of $f$.  In the space spanned by the essential variables, the Laurent semi-norm for polynomials with at least two terms is shown to be a norm.
\end{abstract}

\maketitle

\section{Introduction}
\subsection{Laurent norms}
Given a Laurent polynomial $f$ with integer coefficients in $n$ variables, 
\begin{displaymath} f \in \mathbb{Z}[t_1^{\pm1}, \dots , t_n^{\pm1}], \end{displaymath}
we define a norm $\| \hspace{2mm} \|_f$ on the vector space $(\mathbb{R}^n)^*$. If we write $f$ in multi-index notation, then \begin{displaymath} f = \sum_{\alpha \in \supp(f)} c_{\alpha} \textbf{t}^{\alpha} \end{displaymath}  with $\textbf{t}^{\alpha} = t_1^{\alpha_1} \dots t_n^{\alpha_n}$ and $\supp(f)$ denoting support of  $f$; the set of all $\alpha=(\alpha_1,\dots ,\alpha_n) \in \mathbb{Z}^n$ labeling non-zero constants $c_{\alpha}$. Let $\phi$ be an arbitrary vector in the space dual to the space containing the exponents, the vector space $(\mathbb{R}^n)^*$ dual to $\mathbb{R}^n$.  The vector $\phi$ acts on the integer-valued exponents $\alpha \in \mathbb{Z}^n \subset \mathbb{R}^n$ by duality. The \textit{Laurent norm} of $\phi$ for the Laurent polynomial $f$, denoted $\| \phi \|_f$, is the supremum of $ \phi(\alpha - \beta)$ taken over all the exponents $\alpha, \beta$ in the support of $f$;\begin{displaymath} \| \phi \|_f := \sup_{\alpha,\beta \in \supp f}\phi(\alpha - \beta). \end{displaymath}
In this fashion, every $f \in \mathbb{Z}[t_1^{\pm1}, \dots , t_n^{\pm1}]$ determines a Laurent norm $\| \hspace{2mm} \|_f$ on $(\mathbb{R}^n)^*$. The Laurent norm is sometimes \textit{degenerate}; it can be zero for non-zero vectors $\phi$.

\subsection{The width function for the Newton polyhedron of $f$ }In the theory of polyhedra, there is a standard function called the \textit{width function w}. It has values in the non-negative real numbers and acts on the product space of polyhedra and vectors of the dual space to the ambient space containing these polyhedra.  The \textit{Newton polyhedron $\Ne(f)$} is the convex hull of the set of exponents of $f$.\footnote{The term \textit{Newton polytope} can also be used for the Newton polyhedron.} Let $x$ and $y$ be arbitrary points of the Newton polyhedron in $\mathbb{R}^n$ and let $\phi$ be a vector of the dual vector space $(\mathbb{R}^n)^*$. The width function $w$ for $\Ne(f)$ of the vector $\phi$ is given by 
 \begin{equation}\nonumber w(\Ne(f), \phi): = \sup _{  x,y \in \Ne(f) } \phi(x-y).  
  \end{equation} 
The supremum of $ \phi(x-y)$ must take its values on the vertices of $\N(f)$. This implies the width function of $\Ne(f)$ and $\phi$ is equal to the Laurent norm for $f$ of $\phi$, so that \begin{displaymath} w(\Ne(f), \phi) =\|\phi \|_f. \end{displaymath} 
\subsection{Decomposability of the Laurent norm} By viewing the Laurent norm $\| \cdot \|_f$ for $f$ as the width function $w(\Ne(f),\cdot)$ of the Newton polyhedron $\Ne(f)$ of $f$, a great simplification in the calculation of the Laurent norm can be obtained by applying the theory of \textit{Minkowski sums}.  If $P$ and $Q$ are two polyhedra, the Minkowski sum of $P$ and $Q$ is the polyhedron made up of all vector sums $x + y$ where $x$ is an arbitrary vector of $P$ and $y$ is an arbitrary vector of $Q$.  It is a standard result of the theory of Minkowski sums \cite{Sch} that the width function is a \textit{Minkowski linear} function.  Using this linearity of the width function, we obtain our first main result, the  \textit{Laurent norm decomposition formula}: \begin{theorem} Let  $f$ be a multivariable Laurent polynomial with integer coefficients that can be expressed in terms of its irreducible factors as $f = f_1^{n_1} \cdots f_k^{n_k}$ with $n_i \in \mathbb{N}$, for $i=1, \dots k$. Then the Laurent norm for $f$ of $ \phi$ is given by 
\begin{equation}\nonumber \| \phi \|_{f} = \sum_{i=1}^k n_i \| \phi \|_{f_i}. \end{equation}\end{theorem}

\subsection{ The Laurent norm unit ball} 
We introduce the \textit{essential variables} of the Laurent norm unit ball $\mathcal{B}_f$ for the Laurent polynomial $f$ and use them to determine the reduced Laurent norm unit ball $\tilde{\mathcal{B}}_f$ for $f$. The space spanned by the essential variables is $(\mathbb{R}^m)^*$, where $m$ is the dimension of the Newton polyhedron $\Ne(f)$ and  $m \leq n$. We show that the two unit balls $\mathcal{B}_f$ and $\mathcal{\tilde{B}}_f$ are related by  \begin{equation}\nonumber \mathcal{B}_f =  \mathcal{\tilde{B}}_f \times (\mathbb{R}^{n-m})^*.\end{equation} 
We then deduce our second main result: 
\begin{theorem} The reduced Laurent norm unit ball $\mathcal{\tilde{B}}_f$ is a convex polyhedron. \end{theorem} 
Using this theorem and the coincidence of the Laurent norm for a polynomial $f$ and the width function of the Newton polyhedron $\Ne(f)$ of $f$, we deduce our third main result: 
\begin{theorem} In the space $(\mathbb{R}^m)^*$ spanned by the essential variables, the Laurent norm for a polynomial $f$ with a least two terms is non-degenerate; $ \| \phi\| _f = 0$ if and only if $\phi=0$.\end{theorem}
This theorem implies that in $(\mathbb{R}^m)^*$ with $m \geq 1$, the Laurent norm is not just a semi-norm (degenerate norm) but a norm (non-degenerate norm).

\subsection{Applications of the Laurent norm to the Alexander norm}
The motivation for introducing this new norm comes from a norm defined on the first cohomology group $H^1(M;\mathbb{R})$ of a connected, compact, orientable, 3-manifold $M$ whose boundary (if any) is a union of tori. This norm, introduced by McMullen in \cite{McM}, is directly determined by the multivariable Alexander polynomial of the $3$-manifold $M$, and is called the \textit{Alexander norm} of $M$. The Alexander norm is a special case of the Laurent norm in that the Laurent polynomial $f$ is also the Alexander polynomial $\triangle$ for the 3-manifold $M$. Thus our results above for the Laurent norm can also be applied directly to the Alexander norm. 

We conclude with two examples of calculations of the Alexander norm and unit ball. The first example, which can also be found in McMullen \cite{McM}, is the well-known three component link, the Borromean rings. We find there are no inessential variables and we show, by methods differing from McMullen's in \cite{McM}, that the Alexander norm unit ball is the unit cube. Our second example is original; it is a six-component great circle link. We find the two essential and four inessential variables for the Alexander polynomial of this link and show its reduced Alexander norm unit ball is a unit square.

\section{The Laurent norm}
 Let $f$ be a Laurent polynomial with integer coefficients, \begin{displaymath} f \in  \mathbb{Z}[t_1^{\pm 1}, \dots , t_n^{\pm 1}],\end{displaymath} so that $f$ can be expressed as \begin{eqnarray}\label{E:pol}  f(t_1, \dots ,t_n)&=& \sum c_{\alpha_1 \alpha_2 \dots \alpha_n} t_1^{\alpha_1}t_2^{\alpha_2} \dots t_n^{\alpha_n}\\ \nonumber
&=& \sum_{\alpha \in \supp(f)} c_{\alpha} \textbf{t}^{\alpha} 
 \end{eqnarray}  $\mbox{where } \supp(f) = \lbrace \alpha \colon c_{\alpha} \neq 0 \rbrace$. Then $\alpha \in \mathbb{Z}^n, \forall \alpha \in \supp(f)$.
Further let $\phi \in (\mathbb{R}^n)^* = \Hom_{\mathbb{R}}(\mathbb{R}^n, \mathbb{R}) \cong \mathbb{R}^n$ be an element of the dual vector space $(\mathbb{R}^n)^*$ to the space $\mathbb{R}^n$ in which each $\alpha$ lies in since $\mathbb{Z}^n \subset \mathbb{R}^n$.  \begin{defn}\label{D:fn} With $f$ as above, the \textit{Laurent norm} for $f$ is  \begin{eqnarray}\label{E:fn}  && \| \hspace{2mm} \|_f \colon (\mathbb{R}^n)^* \to \mathbb{R}^+ \cup \{0\}. \\\nonumber
 &&  \| \phi \|_f = \sup_{\alpha,\beta \in \supp(f)} \phi(\alpha - \beta). \end{eqnarray}
\end{defn}
\begin{rem} According to this definition, the Laurent norm for all monomials is identically zero.\end{rem}
\begin{prop} The function $\| \cdot \|_f$ of this definition is a well-defined semi-norm.\end{prop}
\begin{proof}We show that it satisfies the three properties of a semi-norm;\begin{enumerate}
\item $\| \phi \|_f \geq 0$: If $\phi(\alpha - \beta)$ is negative, then $\phi(\beta-\alpha)$ is positive so the supremum which occurs in Equation \eqref{E:fn} is always positive or zero.
\item $\| \phi + \phi' \|_f \leq \| \phi \|_f + \| \phi'\|_f$: Assume that $\| \phi + \phi' \|_f = \phi(\alpha_o - \beta_o)$, then clearly as supremums taken over pairs of exponents of $f$, $\| \phi \|_f \geq \phi(\alpha_o- \beta_o)$ and $\| \phi'\|_f  \geq \phi'(\alpha_o - \beta_o)$.
\item $\| \lambda \phi \|_f = | \lambda | \|\phi \|_f, \forall \lambda \in \mathbb{R}$: For $\lambda \geq 0$ we have that,  \begin{displaymath}\| \lambda \phi \|_f =  \sup_{\alpha,\beta \in \supp(f)} \lambda \phi(\alpha - \beta) = \lambda  \sup_{\alpha,\beta \in \supp(f)} \phi(\alpha - \beta) = \lambda \| \phi \|_f \end{displaymath} and for $\lambda <0$, we find that \begin{displaymath}\| \lambda \phi \|_f =  \sup_{\alpha,\beta \in \supp(f)} \lambda \phi(\alpha - \beta) = -\lambda  \sup_{\alpha,\beta \in \supp(f)} \phi(\alpha - \beta) = -\lambda \| \phi \|_f. \end{displaymath}
\end{enumerate}\end{proof}
For some polynomials $f$, the Laurent norm may be \textit{degenerate}; it takes the value zero for non-zero vectors $\phi$ (See the second example of Section \ref{S:E}.).\footnote{ Since the Laurent norm is sometimes a \textit{degenerate norm}, to be precise, it should be called a \textit{semi-norm}.  However, it is traditional to refer to semi-norms of this type, for example the Alexander and Thurston norms, as norms rather than semi-norms.} We show in Section \ref{S:LB}, that for polynomials with at least two terms, there is a special set of coordinates, which we call \textit{essential variables}, such that the Laurent norm becomes a norm rather than only semi-norm in these coordinates.  
\subsection{The width function for polyhedra and the Laurent norm} 
In Schneider \cite{Sch}, the definitions of two functions for polyhedra, the support function and the width function, can be found. 
Assume that $P$ is a polyhedron in $\mathbb{R}^n$, and that as before $\phi$ is a vector of the dual vector space $(\mathbb{R}^n)^*$.  The support and width functions, $h$ and $w$ respectively, of $P$ and the vector $\phi$ are given by the following.\footnote{In the definition of the support and width functions for polyhedra (for example in \cite{Sch}), it is traditional to assume that the vector $\phi$ on which they are defined must be a vector of length one for the standard Euclidean norm. We do not impose this restriction. Also in the theory of polyhedra, the duality relation between the vectors $\phi$ and $x$, $\phi(x)$, is just the standard Euclidean scalar product between $\phi$ and $x$ in $\mathbb{R}^n$.}
 \begin{defn}\label{D:SF}The \textit{support function $h$} of the polyhedron $P$ and the vector $\phi$ is 
\begin{equation}\nonumber 
h(P,\phi) = \sup_{  x \in P}   \phi(x) .   \end{equation} \end{defn}
\begin{defn}\label{D:WF}The \textit{width function $w$} of the polyhedron $P$ and the vector $\phi$ is defined as 
\begin{equation} w(P, \phi) = \sup _{  x,y \in P } \phi(x-y).  
  \end{equation} \end{defn}

\subsubsection{Geometric interpretation of the width function } If $x$ and $y$ lie on the boundary of the polyhedron $P$, then $z = x-y$ is a ``width'' vector, with length equal to the width of $P$ along the direction determined by the vector $z$. Then $\phi(z)$ is the projection of this width vector $z$ along the vector $\phi$.  Since $w(P,\phi)$ is the supremum of all the widths of $P$ along the direction determined by $\phi$ up to a scale factor, it is called the width function for $P$. 

Next, we show that the width function of the Newton polyhedron $\Ne(f)$ of Laurent polynomial $f$ is equal to the Laurent norm for $f$.
\begin{prop}\label{P:W} The Laurent norm for a Laurent polynomial $f$ is equal to the width function of the Newton polyhedron $\Ne(f)$, so that \begin{displaymath} \| \phi \|_f = w(\Ne(f), \phi). \end{displaymath}\end{prop}

\begin{proof} The proof follows directly from Definition \ref{D:fn} of the Laurent norm for $f$ and Definition \ref{D:WF} of the width function of $\Ne(f)$. These two definitions are the same except in the definition of the width function the supremum is taken over all the points of the convex hull that make up the polyhedron $\N(f)$, while in the definition of the Laurent norm the supremum is taken only over a finite generating set of points containing all the vertices of the convex hull of the Newton polyhedron $N(f)$.  However, since the supremum which occurs in the width function always takes its values on the vertices of the convex polyhedron $\Ne(f)$, the width function $w(\Ne(f),\phi)$ of $\Ne(f)$ must have the same value as the Laurent norm $\| \phi \|_f$ for $f$ for all $\phi$. \end{proof}
\subsection{The single variable polynomial method of calculating the Laurent norm}
In \cite{McM}, McMullen states that the Alexander norm  for the multivariable Alexander polynomial $\triangle(t_1, \dots ,t_r)$ can be easily calculated using the single variable Alexander polynomial, which he denotes $\triangle^{\phi}$. 
We can easily generalize this method to calculation of the Laurent norm.
\begin{defn}\label{D:McAP} Assume we are given a multivariable Laurent polynomial $f(t_1, \dots , t_n)$. We define an associated \textit{single variable polynomial $f^{\phi}(t)$ }that depends also on a dual vector $\phi$ of $(\mathbb{R}^n)^*$, as follows: \begin{equation} \nonumber f^{\phi}(t) =  f(t^{\phi_1}, \dots ,t^{\phi_n}).\end{equation} \end{defn}

In terms of this single variable polynomial $f^{\phi}(t)$, the Laurent norm of $\phi$ for $f$ is given simply by \begin{equation}\nonumber \| \phi \|_f = \deg f^{\phi}(t), \end{equation} where $\deg$ denotes the degree of the single variable polynomial, which is the difference between the largest and smallest exponent of all the terms in the polynomial. 
We show that this expression for the Laurent norm in terms of the degree of the single variable polynomial $f^{\phi}$ is equivalent to original Definition \ref{D:fn} in the next proposition. 
\begin{prop}\label{P:Mman} Let $f(t_1, \dots ,t_n)$ be a multivariable Laurent polynomial, and assume that $f^{\phi}(t) \neq 0$ for the vector $\phi \in (\mathbb{R}^n)^*$. Then the Laurent norm of $\phi$ for $f$ is given by $\| \phi \|_f = \deg f^{\phi}(t)$.\end{prop}
\begin{proof} The vector $\phi$ induces a map, which we denote as $\phi_*$, from the ring of multivariable Laurent polynomials $\mathbb{Z}[t_1^{\pm1}, \dots ,t_n^{\pm 1}]$ to the ring of single variable Laurent polynomials $\mathbb{Z}[t^{\pm 1}]$: \begin{eqnarray} \nonumber  && \phi_* \colon \mathbb{Z}[t_1^{\pm1}, \dots ,t_n^{\pm 1}] \to \mathbb{Z}[t^{\pm 1}], \\ \nonumber
&& \phi_*(t_i) = t^{\phi_i}, i=1, \dots , n.  \end{eqnarray}
This map acts on the multivariable  polynomial $f(t_1, \dots ,t_n)$ to give the single variable polynomial $f^{\phi}(t)$ for the vector $\phi$:\begin{equation} \nonumber f^{\phi}(t)= \phi_*(f(t_1, \dots ,t_n))=f(t^{\phi_1}, \dots ,t^{\phi_n}).\end{equation} If write the polynomial $f$ in the multi-index notation as \begin{equation} \nonumber f(t_1, \dots ,t_n) = \sum_{\alpha \in \supp(f)} c_{\alpha}\textbf{t}^{\alpha}, \end{equation} then the Laurent norm of $\phi$ for $f$ can be expressed using $\phi_*$ as  \begin{equation}\label{E:Anorm} \| \phi \|_f = \sup_{\alpha,\beta \in \supp(f)} | \phi_*(\textbf{t}^{\alpha} - \textbf{t}^\beta)| . \end{equation} In this expression of the Laurent norm, $\phi_*$ acts on the generator $\textbf{t}^{\alpha } = t_1^{\alpha_1} \cdots t_n^{\alpha_n}$ instead of the exponent $\alpha$ of the generator as is the case for $\phi$. Since \begin{eqnarray} \nonumber   \phi_*(\textbf{t}^{\alpha})&=& (t^{\phi_1})^{\alpha_1} \cdots (t^{\phi_n})^{\alpha_n}\\ \nonumber & =& t^{\phi_1 \cdot \alpha_1 + \cdots +\phi_n \cdot \alpha_n}, \end{eqnarray} the supremum in Equation \eqref{E:Anorm} is the degree of the single variable polynomial $f^{\phi}(t)$, as claimed. \end{proof}
\begin{rem} This method for calculating the Laurent norm is only valid for the vectors $\phi$ for which $f^{\phi}(t) \neq 0$. By this method, if $f^{\phi}(t)=0$, then  necessarily $\| \phi \|_f = 0 $. However, when the original Definition \ref{D:fn} of the Laurent norm for $f$ is applied instead, we find that in some instances we obtain the contradictory result that $\| \phi \|_f \neq 0$ even though $f^\phi(t) = 0$.  \end{rem}

\section{The Laurent norm decomposition formula}
In this section, we use the notions of Minkowski sum and Minkowski linearity to derive an expression for the Laurent norm of the polynomial $f$ in terms of the Laurent norms of each of its irreducible factors.

A very extensive theory of Minkowski sums can be found in Schneider \cite{Sch}. 
\begin{defn}[Gelfand et al. \cite{Gel}, Def.~ 5.4.7] The \textit{Minkowski sum} of two polyhedra $P$ and $Q$ is the set of all vector sums $x + y$ with $x \in P$ and $y  \in Q$.\end{defn}

We shall make use of the following two standard results of the theory of polyhedra.

\begin{prop}[Gelfand et al. \cite{Gel}, Prop.~6.1.2(b)]\label{P:MS} The Newton polyhedron, $\Ne(fg)$, of the product of two polynomials $f$ and $g$ equals the Minkowski sum of the Newton polyhedra of the factors, $\Ne(f) + \Ne(g)$:
\begin{equation}\nonumber \Ne(fg) = \Ne(f) + \Ne(g). \end{equation}\end{prop}
By applying this proposition inductively to the product of $f$ with itself we obtain the following corollary.
\begin{cor}\label{CO:N} The Newton polyhedron $\Ne(f^n)$ of a polynomial $f$ to an integer power $n \in \mathbb{N}$, $f^n$, is the product of the power by the Newton polyhedron of the polynomial, $n \cdot \Ne(f):$
\begin{equation}\nonumber \Ne(f^n) = n \cdot \Ne(f).\end{equation}\end{cor}

 Since Minkowski linearity is not the same as the usual form of linearity we define it explicitly.
\begin{defn}A real-valued function $g$ from polyhedra to the real numbers is said to be \textit{Minkowski linear} if it satisfies the following two properties: \begin{enumerate}
\item{Minkowski additivity: }$g$ is Minkowski additive if it satisfies $g(Q+P) = g(Q) + g(P)$, where $Q$ and $P$ denote arbitrary polyhedra. 
\item{Minkowski scaling: }$g$ satisfies the Minkowski scaling property if $g(\lambda \cdot P) = \lambda g(P)$, for any $\lambda \in \mathbb{R}^+$.
\end{enumerate}
\end{defn}
It is well-known that the width function $w$ is Minkowski linear on polyhedra. 
\begin{theorem}\label{L:Min} The width function is a Minkowski linear function on polyhedra. That is it satisfies the following two properties. \begin{enumerate} \item{ Minkowski additivity:} $ w(P + Q,\phi) = w (P, \phi)  + w(Q,\phi)$.
\item{Minkowski scaling: }$w(\lambda \cdot P,\phi)  = \lambda \cdot w(P,\phi),   \forall \lambda \in \mathbb{R}^+ $.\end{enumerate}
\end{theorem}
A detailed proof of this theorem can be found in Long \cite{Lont}. A proof that the support function is a Minkowski additive function can be found in Theorem 1.75 of \cite{Sch}.  It is then stated in \cite{Sch} that Minkowski additivity of the support function implies that the width function is also a Minkowski additive function. 

Since the Laurent norm for $f$ is equal to the width function for $\Ne(f)$, by applying the above theorem to the Laurent norm, we obtain the following corollary. 
 
\begin{cor}\label{C:ML} The Laurent norm for $f$ satisfies the following properties:\begin{enumerate}
\item  $\| \phi \|_{f\cdot f'} = \| \phi \|_f + \| \phi \|_{f'}$.
\item  $\| \phi \|_{f^n} = n \| \phi \|_f, \forall n \in \mathbb{N}.$
\end{enumerate}
\end{cor}
\begin{proof}  For any Laurent polynomial $f$, $\| \phi \|_{f} = w(\Ne(f), \phi)$ by Proposition \ref{P:W}.  By Proposition \ref{P:MS}, $\Ne(f\cdot f') = \Ne(f) + \Ne(f')$ .  Applying the Minkowski additivity property of $w$, we obtain \begin{eqnarray} \nonumber   \| \phi \|_{f\cdot f'} & = & w(\Ne(f \cdot f'), \phi)\\ \nonumber & =& w(\Ne(f) + \Ne(f'), \phi)\\ \nonumber &
=& w(\Ne(f),\phi)+ w(\Ne(f'), \phi) \\ \nonumber &
=& \| \phi \|_f + \| \phi \|_{f'}. \end{eqnarray}
 By Corollary \ref{CO:N}, $\Ne(f^n) = n \cdot \Ne(f)$. Applying the Minkowski scaling property of $w$, we obtain
\begin{eqnarray} \nonumber  \| \phi \|_{f^n} &
= & w(\Ne(f^n), \phi) \\ \nonumber &
= & w(n \Ne(f), \phi)\\ \nonumber &
= & n \cdot w(\Ne(f), \phi)\\ \nonumber &
= & n \| \phi \|_f.  \end{eqnarray}
\end{proof}
Applying this corollary inductively to an arbitrary Laurent polynomial $f$ expressed as a product of its irreducible factors $f_1, \dots , f_k$, we obtain our main result, \textit{the Laurent norm decomposition formula}. 
\begin{theorem} Let the multivariable Laurent polynomial $f$ with integer coefficients be written in terms of its irreducible factors, $f_1, \dots , f_k$, as  $f = f_1^{n_1} \cdots f_k^{n_k}$ with $n_i \in \mathbb{N}, i=1, \dots k$. Then the Laurent norm of $\phi$ for this polynomial is given in terms of the Laurent norms of $\phi$ for the irreducible factors of $f$ by 
\begin{equation}\label{E:f1} \| \phi \|_{f} = \sum_{i=1}^k n_i \| \phi \|_{f_i}. \end{equation}\end{theorem}
\begin{proof}   If $k=1$ and $f=f_1^{n_1}$, then \begin{displaymath}\| \phi \|_{f} = n_1 \| \phi \|_{f_1}\end{displaymath} by the second property of Corollary \ref{C:ML}. Hence Equation \eqref{E:f1} is true for $k=1$. Assume that Equation \eqref{E:f1} for $f$ is true for $k=m-1$ and $k>1$. We show this implies it is also true for $k=m$. Set $\tilde{f}= f_1^{n_1} \cdots f_{m-1}^{n_{m-1}}$ so that $f = \tilde{f} \cdot f_m^{n_m}$. The first property of Corollary \ref{C:ML} implies \begin{equation} \label{E:f3} \| \phi \|_f = \| \phi \|_{\tilde{f}} + \| \phi \|_{f_m^{n_m}}.\end{equation}  By the inductive hypothesis, the Laurent norm $\| \phi \|_{\tilde{f}}$ of polynomial $\tilde{f}$  satisfies Equation \eqref{E:f1} with $k=m-1$. By the second property of Corollary \ref{C:ML},  \begin{displaymath} \| \phi \|_{f_k^{n_m}}= n_m \| \phi \|_{f_m} \end{displaymath} After substituting this result into Equation \eqref{E:f3}, we find \begin{equation}\nonumber \label{E:f2} \| \phi \|_f = \| \phi \|_{\tilde{f}} + n_m \| \phi \|_{f_m}. \end{equation}We conclude by induction that Equation \eqref{E:f1} is true for every $k$. \end{proof} 
\section{Two forms of duality and essential variables}
A polyhedron which has dimension $m < n$ in the space $\mathbb{R}^n$ has a dual in $(\mathbb{R}^n)^*$, which is not the same as its dual in the space of the same dimension, $\mathbb{R}^m$. For this reason, we introduce two types of duality for a polyhedron $P$; one in terms of the ambient space $\mathbb{R}^n$ that $P$ is in and the other in terms of the subspace $\mathbb{R}^m$ of the ambient space $\mathbb{R}^n$ which is of the same dimension as $P$. 
\begin{defn}\label{D:DPo} If $P$ is a polyhedron in $\mathbb{R}^n$ such that the origin is in its interior, the \textit{dual polyhedron $P^*$} is defined as  \begin{equation}\nonumber P^* = \lbrace \phi  \in  (\mathbb{R}^n)^* \mid \phi(x) \leq 1, \forall  x  \in  P \rbrace. \end{equation} \end{defn}

This definition can be applied to the Newton polyhedron $\Ne(f)$ of the polynomial $f$ to determine the \textit{dual Newton polyhedron $\Ne(f)^*$}.

\subsection{Essential variables} Assume the polynomial $f$ has $n$ variables, but its Newton polyhedron has dimension $1 \leq m \leq n$.  Then there is a coordinate system for which each point of the $m$-dimensional Newton polyhedron $\Ne(f)$ has $m$ non-zero coordinates and $n-m$ coordinates identically zero.  We call these $m$ coordinates the \textit{essential variables} of the Newton polyhedron. Without loss of generality, we can assume that the first $m$ coordinates of $\mathbb{R}^n$ are the essential variables, $(x_1, \dots, x_m)$. These determine $m$ essential variables $(\phi_1, \dots,\phi_m)$ of the dual space $(\mathbb{R}^n)^*$ . 
\begin{prop}\label{P:E1} The Laurent norm $\| \phi \|_f$ of $\phi \in (\mathbb{R}^n)^*$ for the $n$-variable polynomial $f$ is a function of only the $m$ essential variables, $(\phi_1, \dots ,\phi_m) \in (\mathbb{R}^m)^*$.\end{prop}\begin{proof} In the coordinate system for $\mathbb{R}^n$ of the $m$ essential variables described above, the polynomial $f$ is a function of only $m$ variables $(t_1, \dots , t_m)$. Each of its exponents is a function of the $m$ essential variables $(x_1, \dots , x_m) \in \mathbb{R}^m$. It follows that the Laurent norm for $f$ which is directly determined by the exponents of $f$ must be a function of only the $m$ essential variables $(\phi_1, \dots, \phi_m) \in (\mathbb{R}^m)^*$.\end{proof} 
We can define a second duality transformation of the Newton polyhedron which is the duality transformation restricted to the $m$ essential variables.  We call this the \textit{reduced dual Newton polyhedron}, and denote it by  $\tilde{\Ne}(f)^*$.
\begin{defn}\label{D:Po1} \textit{The reduced dual $\tilde{\Ne}(f)^*$ to the Newton polyhedron $\Ne(f)$} is determined by the equation
\begin{equation}\nonumber  \tilde{\Ne}(f)^* = \lbrace \phi  \in  (\mathbb{R}^m)^* \mid \phi(x) \leq 1, \forall x   \in \Ne(f) \rbrace. \end{equation} It is the dual of the Newton polyhedron in the space of essential variables.\end{defn}

\begin{prop}\label{P:Du1} The reduced dual Newton polyhedron $\tilde{\Ne}(f)^*$ is a convex $m$-dimensional polyhedron.\end{prop}
\begin{proof} When we restrict the duality transformation to the space of essential variables $\mathbb{R}^{m}$, the dual $\tilde{\Ne}(f)^*$ to the Newton polyhedron $\Ne(f)$, which is convex by definition, must be a convex polyhedron of the same dimension $m$.  This is a well-known property of the duality transformation for polyhedra.\end{proof}
\begin{prop}\label{P:Du2} The dual Newton polyhedron $\Ne(f)^*$ is related to the reduced dual Newton polyhedron $\tilde{\Ne}(f)^*$  by \begin{equation}\nonumber \label{E:Du1} \Ne(f)^* =  \tilde{\Ne}(f)^* \times (\mathbb{R}^{n-m})^*. \end{equation}\end{prop}
\begin{proof} The dual Newton polyhedron $\Ne(f)^*$ is the set of all $\phi \in (\mathbb{R}^{n})^*$ which satisfy the condition of the inequality of the duality transformation in $\mathbb{R}^n$, $\phi(x) \leq 1$ for all $x \in \Ne(f)$. The set depends only on the values of the essential coordinates because the inessential coordinates $(x_{m+1},  \dots ,x_{n})$ of $\Ne(f)$ are identically zero. Hence for all $x \in \Ne(f)$, we have that $\phi_i x_i = 0, i=m+1, \dots ,n $ for every  $\phi_i \in \mathbb{R}^*$.  This implies each inessential variable contributes a factor of $\mathbb{R}^*$ to the dual to the Newton polyhedron as a set in $(\mathbb{R}^{n})^*$ as claimed. Further, by the above Proposition \ref{P:Du1}, the set of essential variables contribute the factor $\tilde{\Ne}(f)^*$. Hence $\Ne(f)$ is the product set of these two sets, $\tilde{\Ne}(f)^*$ and $(\mathbb{R}^{n-m})^*$.  \end{proof}
\section{The Laurent norm unit ball} \label{S:LB} 
In this section we describe some properties of the unit ball in the Laurent norm.
\begin{defn}The \textit{Laurent norm unit ball $\mathcal{B}_f$} for polynomial $f$ is the set \begin{equation}\nonumber \mathcal{B}_f = \lbrace \phi \in (\mathbb{R}^n)^* \mid \| \phi \|_f \leq 1 \rbrace. \end{equation}\end{defn}
In the coordinate system containing the essential coordinates, the reduced Laurent norm unit ball for $f$ can be defined as follows.
\begin{defn} The \textit{reduced Laurent norm unit ball  for $f$} is the set \begin{equation}\nonumber \tilde{\mathcal{B}}_f = \lbrace \phi \in (\mathbb{R}^m)^* \mid \| \phi \|_f \leq 1 \rbrace, \end{equation} where $m$ is the dimension of the Newton polyhedron, and $(\mathbb{R}^m)^*$ is the space spanned by the essential coordinates for $\tilde{\Ne}(f)^*$.\end{defn}
Using these two definitions, we are able to directly relate the two unit balls $\mathcal{B}_f$ and $\tilde{\mathcal{B}}_f$ in the following theorem.
\begin{theorem} \label{T:Rln}The Laurent norm unit ball $\mathcal{B}_f$ is related to the reduced Laurent norm unit ball $\tilde{\mathcal{B}}_f$ by \begin{equation} \nonumber \mathcal{B}_f = \tilde{\mathcal{B}}_f \times (\mathbb{R}^{n-m})^*.  \end{equation}\end{theorem}
\begin{proof} In the coordinate system containing the $m$ essential variables, the Laurent norm for $f$ of a vector $\phi \in (\mathbb{R}^n)^*$ is a function of only the $m$ essential coordinates by Proposition \ref{P:E1}. Hence the inequality $\| \phi \|_f \leq 1$ does not depend on the value of any of the inessential coordinates and each of them contributes a factor of $\mathbb{R}^*$ to $\mathcal{B}_f$. Since there are $n-m$ inessential coordinates, their total contribution is $(\mathbb{R}^{n-m})^*$. In the space spanned by the $m$ essential variables, $(\mathbb{R}^m)^*$, the set of points which satisfy $\| \phi \|_f \leq 1$ is given by the reduced Laurent norm unit ball $\mathcal{\tilde{B}}_f$ by construction. We conclude that $\mathcal{B}_f$ is the product set of $\mathcal{\tilde{B}}_f$ and $(\mathbb{R}^{n-m})^*$. 
 \end{proof}

We now prove that the reduced Laurent norm unit ball is a convex polyhedron using the following lemma and some standard facts about polyhedra. 
\begin{lem}\label{C:Rub}The reduced Laurent norm unit ball $\tilde{\mathcal{B}}_f$ for $f$ is a bounded set, for the Euclidean norm in $(\mathbb{R}^m)^*$, of dimension $m$.\end{lem}
\begin{proof}  The Laurent norm $\| \cdot \|_f$ is the width function $w(\Ne(f), \cdot)$ of the $m$-dimensional Newton polyhedron $\Ne(f)$ by Proposition \ref{P:W}. Hence, it is geometrically obvious that the reduced Laurent norm unit ball $\tilde{\mathcal{B}}_f$ must be a bounded set in $(\mathbb{R}^m)^*$ of dimension $m$ for the following reasons.  \begin{enumerate} \item There is a finite-length (for the Euclidean norm in $\mathbb{R}^m$) width vector $z=x-y$ with $x,y \in \Ne(f)$ of the polyhedron $\Ne(f)$ in every direction. \item The projection $\phi(z)$ of any finite-length (for the Euclidean norm in $(\mathbb{R}^m)^*$) dual vector $\phi$ onto the vector $z$ must be also be finite. \item The width function $w(\Ne(f), \phi)$, which is the supremum of all the projections $\phi(z)$ for the fixed  vector $\phi$, is finite and non-zero. Let us say that $w(\Ne(f), \phi) = c$ where $c$ is a positive finite real number.
\item The boundary of the reduced Laurent norm unit ball in the direction of $\phi$ is given by the dual vector $\tilde{\phi} = c^{-1} \cdot \phi$ and this vector $\tilde{\phi}$ is of finite Euclidean length because $\phi$ is of finite length and $c$ is finite.
\item Since the boundary of the reduced norm unit ball consists entirely of vectors of finite Euclidean length, it is bounded as a subset of $(\mathbb{R}^m)^*$.
\item The unit ball must be of dimension $m$. If we choose a coordinate system for which the $m$-dimensional Newton polyhedron $\Ne(f)$ contains the origin, then it is obvious that the width function must be non-zero along every direction in $(\mathbb{R}^m)^*$.\end{enumerate} \end{proof}

By definition, a convex polyhedron is the convex hull of a finite set of points. A convex $m$-dimensional polyhedron $P$ can be also characterized as the intersection of a set of half-spaces $H_{\alpha}$ in $\mathbb{R}^m$, one for each of the top-dimensional faces called facets $\mathcal{F}_{\alpha}$ of $P$. 
The boundary $\partial H_{\alpha}$ of each half-space $H_{\alpha}$ is an $(m-1)$-dimensional hyperplane $\mathcal{H}_{\alpha}$. Since $\mathcal{F}_{\alpha} \subset \mathcal{H}_{\alpha}$, the hyperplane $\mathcal{H}_{\alpha}$ is called the \textit{supporting hyperplane} to the facet $\mathcal{F}_{\alpha}$.

Using these facts and Lemma \ref{C:Rub}, we obtain the following theorem:
\begin{theorem}\label{T:Cp} The reduced Laurent norm unit ball $\tilde{\mathcal{B}}_f$ is a convex polyhedron. \end{theorem}
\begin{proof} 
All vectors $\phi$ along the same ray $\mathcal{R}$ from the origin in $(\mathbb{R}^m)^*$  are equal up to a scale factor and this equivalence implies their Laurent norms are also equal up to the same scale factor; if $\tilde{\phi}= c \phi$ with $c \in \mathbb{R}^+$, then $ \| \tilde{\phi}\|_f = c \| \phi \|_f$. This implies that the supremum which occurs in the Laurent norm takes its value, $\phi(\alpha_o - \beta_o)$, on the same pair of vertices $(\alpha_o,\beta_o)$ of $\tilde{\Ne}(f)$ for all $\phi \in \mathcal{R}$.  Since the set $\textbf{F}$ of all the pairs of vertices on which the Laurent norm takes its value is finite, this supremum must also occur on the same pair of vertices $(\alpha_o,\beta_o)$ for all the vectors $\phi$ in some open cone $\mathcal{C}_{\alpha_o-\beta_o}$ pointed at the origin. Thus for all points $\phi$ in the cone $\mathcal{C}_{\alpha_o-\beta_o}$, we have that Laurent norm satisfies $\| \phi \|_f=\phi(\alpha_o-\beta_o)$. 

Each pair of vertices $(\alpha,\beta)$ in the set $\textbf{F}$ determines a cone $\mathcal{C}_{\alpha-\beta}$. The cone $\mathcal{C}_{\alpha-\beta}$ intersects the $(m-1)$-dimensional hyperplane 
\begin{equation} \nonumber \mathcal{H}_{\alpha-\beta} = \{ \phi \in (\mathbb{R}^m)^* \mid  \phi(\alpha - \beta) = 1\} \end{equation} 
in a set \begin{equation} \nonumber  \mathcal{F}_{\alpha - \beta}= \mathcal{H}_{\alpha-\beta} \cap \mathcal{C}_{\alpha-\beta} \end{equation} on the boundary $\partial \tilde{\mathcal{B}}_f$ of $\tilde{\mathcal{B}}_f$.  

We claim the sets $\mathcal{F}_{\alpha - \beta}$, with $(\alpha,\beta) \in \textbf{F}$ are the facets of a polyhedron.  The hyperplane $\mathcal{H}_{\alpha-\beta}$ determines the half-space $H_{\alpha-\beta}$ in $(\mathbb{R}^m)^*$ which is the set \begin{equation} \nonumber H_{\alpha-\beta} = \{ \phi \in (\mathbb{R}^m)^* \mid  \phi(\alpha - \beta) \leq 1\}. \end{equation} Every point $\phi$ of $\mathcal{\tilde{B}}_f$ must lie in each of these half-spaces since $\| \phi \|_f \leq 1$ for all $\phi \in \mathcal{\tilde{B}}_f$ by definition of $\mathcal{\tilde{B}}_f$, and because as a supremum, $\| \phi \|_f \geq \phi(\alpha-\beta)$ for all $\phi \in (\mathbb{R}^m)^*$ and all pairs $(\alpha,\beta) \in \textbf{F}$. Hence the reduced Laurent norm unit ball is the set \begin{equation} \label{E:Rnub} \mathcal{\tilde{B}}_f = \bigcap_{(\alpha,\beta) \in \textbf{F}} H_{\alpha-\beta}. \end{equation}
Since $\mathcal{\tilde{B}}_f$ is bounded as a set in $(\mathbb{R}^m)^*$ by the above Lemma \ref{C:Rub}, this intersection of half-spaces must be a convex polyhedron with an $(m-1)$-dimensional facet $\mathcal{F}_{\alpha-\beta}$ for each pair $(\alpha,\beta) \in \textbf{F}$ as claimed. \end{proof}
The proof of this theorem, directly implies the following corollary.
\begin{cor} Equation \eqref{E:Rnub} can be used to characterize the reduced Laurent norm unit ball $\mathcal{\tilde{B}}_f$ as an intersection taken over a finite set of half-spaces $H_{\alpha-\beta}$. \end{cor} 
We now use the property that the Laurent norm is the width function of the Newton polyhedron to show that in the space of essential variables the Laurent semi-norm becomes a norm.
\begin{theorem}\label{T:Nd} In the space $(\mathbb{R}^m)^*$, spanned by the essential variables, the Laurent norm for a polynomial $f$ with at least two terms is non-degenerate; $ \| \phi\| _f = 0$ if and only if $\phi=0$.\end{theorem}
\begin{proof} Let us assume that the Newton polyhedron of dimension $m \geq 1$ contains the origin in $\mathbb{R}^m$. Then it is geometrically obvious that the width function $w$ of this polyhedron must be non-zero in every direction, $w(\Ne(f),\phi)\neq 0, \mbox{ for all non-zero } \phi \in (\mathbb{R}^m)^*$, because the $m$-dimensional Newton polyhedron must have a non-zero width as a polyhedron in $\mathbb{R}^m$ in every direction.  Since $w(\Ne(f),\phi)=\| \phi \|_f$ by Proposition \ref{P:W}, the theorem is proved. \end{proof} 
\begin{rem} For a polynomial $f$ which is a monomial, the Newton polyhedron $\Ne(f)$ is zero dimensional, $m=0$, so that the number of essential variables is zero.  Since the Laurent norm $\| \cdot \|_f$ for the monomial $f$ is identically zero, it can not be transformed from a semi-norm to a norm by the correct choice of coordinates.  \end{rem}

\subsection{Calculation of the reduced Laurent norm unit ball}
The following technical lemma make use of the polyhedral property of the reduced Laurent norm unit ball to simplify calculations of the set of its vertices (and other faces of codimension greater than one).

As before, let a Laurent polynomial $f$ be given in terms of its irreducible factors by $f = f_1^{n_1} \dots f_k^{n_k}$. Its Laurent norm on $\phi$ is given by $\| \phi \|_f = \sum_{i=1}^k n_i \| \phi \|_{f_i}$ in essential variables.
\begin{lem}\label{L:CL}  With notation as above, the faces of codimension greater than one of the reduced Laurent norm unit ball $\mathcal{\tilde{B}}_f$ can be determined by setting the Laurent norms for the irreducible components $\| \phi \|_{f_i},$ with $i \in \lbrace 1, \dots ,k \rbrace$, equal to zero for subsets of the set of integers $\lbrace 1, \dots ,k \rbrace$.\end{lem} \begin{proof} The Laurent norms for the irreducible factors, $\| \phi \|_{f_i}, i=1, \dots ,k$, are absolute values of linear functions of the components of $\phi$. The boundary $ \partial \tilde{\mathcal{B}}_f$ of $\mathcal{\tilde{B}}_f$ is the set $\partial \mathcal{\tilde{B}}_f = \lbrace \phi \in (\mathbb{R}^m)^* \mid \| \phi \|_f = 1 \rbrace$.  Each face of codimension greater than one on the boundary of a polyhedron forms a discontinuity of the surface of the polyhedron created by the intersection of the higher dimensional faces.   These discontinuities of the surface of the polyhedron $\tilde{\mathcal{B}}_f$ are formed due to discontinuities of the Laurent norm as a function on the boundary defined by $\| \phi \|_f=1$.   The Laurent norm $\| \phi \|_f$ as a sum of Laurent norms $\| \phi \|_{f_i}, i=1, \dots ,k$, each of which is an absolute value function, is discontinuous at all the points $\phi \in \partial \tilde{\mathcal{B}}_f$ where at least one of these absolute value functions takes the value zero. This fact is a well-know property of the absolute value function.\end{proof}

\section{The Laurent norm for symmetric polynomials}
In general, there is no simple relation between the Laurent norm unit ball $\mathcal{B}_f$ and the dual to the Newton polyhedron $\Ne(f)^*$. However, in the case that the polynomial $f$ satisfies a symmetry invariance under a change of signs of the exponents, the two are related by a scale factor of two as will be proved in the following proposition.
This symmetry property is especially important because the Alexander polynomials are symmetric polynomials.
\begin{defn} The polynomial $f$ is a \textit{symmetric polynomial} if it has the property that $\alpha \in \supp(f)$ if and only if $-\alpha \in \supp(f)$. \end{defn}
In the rest of this section, we assume that $f$ is a symmetric polynomial and that the Newton polyhedron $\Ne(f)$ has been centered at the origin by the appropriate coordinate transformation.
\begin{prop}\label{P:Sc} With $f$ and $\Ne(f)$ as above, the Laurent norm unit ball $\mathcal{B}_f$ for $f$ and the dual polyhedron $\Ne(f)^*$ of the Newton polyhedron $\Ne(f)$ are related by a scale factor of two;\begin{equation}\nonumber  \Ne(f)^* = 2 \mathcal{B}_f.\end{equation}\end{prop}
\begin{proof} The symmetry invariance under a change of sign of the exponents of the polynomial $f$ implies that its Newton polyhedron is symmetric about some central point. Hence if we use a coordinate system in $\mathbb{R}^n$ which locates the central point at the origin, we find that the width function in any direction is just twice the value of the support function in that direction. Since the Laurent norm for $f$ is the width function for the Newton polyhedron $\Ne(f)$, and the support function is determined by the action of duality in the vector space $(\mathbb{R}^n)^*$, we can deduce the following equalities by applying Definitions \ref{D:DPo}, \ref{D:SF} and \ref{D:WF} and Proposition \ref{P:W}. 
\begin{eqnarray}\nonumber  \Ne(f)^* & = & \lbrace \phi  \in  (\mathbb{R}^n)^* \mid \phi(x) \leq 1, \forall x \in  \Ne(f) \rbrace \\ \nonumber 
& = &\lbrace \phi  \in  (\mathbb{R}^n)^* \mid h(\Ne(f),\phi) \leq 1  \rbrace \\ \nonumber
& =& \lbrace \phi  \in  (\mathbb{R}^n)^* \mid w(\Ne(f),\phi) \leq 2  \rbrace\\ \nonumber
& =& \lbrace \phi  \in  (\mathbb{R}^n)^* \mid \|\phi\|_f \leq 2  \rbrace \\ \nonumber
& = &2 \mathcal{B}_f.  \end{eqnarray} \end{proof}
From this proposition, we directly obtain the following corollary.
\begin{cor} \label{C:Sc1} With $f$ and $\Ne(f)$ as above,  the reduced dual Newton polyhedron $\tilde{\Ne}(f)^*$ and the reduced Laurent norm unit ball $\tilde{\mathcal{B}}_f$ for $f$ are related by \begin{equation} \nonumber  \tilde{\Ne}(f)^* = 2 \tilde{\mathcal{B}_f} \end{equation}\end{cor}
\begin{proof} The proof is exactly the same as in the proposition except we use that the vector $\phi$ is in the coordinate system of essential coordinates, $\phi \in (\mathbb{R}^m)^*$, instead of $(\mathbb{R}^n)^*$.\end{proof}
Finally, we find that each vertex $\alpha$ of the Newton polyhedron determines a simple expression for the Laurent norm in the case that the polynomial $f$ is symmetric. 
\begin{prop} Let $f$ and $\Ne(f)$ be as above. Then in essential coordinates the Laurent norm for $f$ of each vector $\phi$ in the cone $\mathcal{C}_{\alpha}$ pointed at the origin through the facet $\mathcal{F}_{\alpha}$ of the reduced Laurent norm unit ball $\tilde{\mathcal{B}}_f$ is given by \begin{equation} \nonumber \| \phi \|_f = |2 \phi(\alpha)|. \end{equation}   \end{prop}
\begin{proof}  By the Definition \ref{D:Po1} of the dual Newton polyhedron $\tilde{\Ne}(f)^*$, a point $\phi'$ on a facet of $\tilde{\Ne}(f)^*$  satisfies $\phi'(\alpha) =1$ because each facet lies on the boundary. On the other hand, a point $\phi$ on the facet $\mathcal{F}_{\alpha}$  of the reduced norm unit ball $\tilde{\mathcal{B}}_f$, along the same ray through the origin as the point $\phi'$, satisfies  $\| \phi \|_f=1$ by definition of the unit ball. By Corollary \ref{C:Sc1}, $\tilde{\Ne}(f)^*=2\tilde{\mathcal{B}}_f$,  so that $\phi' = 2 \phi$. This implies we must have that $\| \phi \|_f = |2 \phi(\alpha)|$ as claimed.  If $\| \phi \|_f=|2 \phi(\alpha)|$ for a point on the boundary of the unit ball, the same result for the Laurent norm applies for all the points on the same ray through the origin because they are equivalent to the boundary point up to a scale factor. We conclude that since our arguments above apply to all points in the facet $\mathcal{F}_{\alpha}$ of $\tilde{\mathcal{B}}_f$, they apply equally well to all the points in the cone $\mathcal{C}_{\alpha}$ through that facet. \end{proof}  

As a direct result of this proposition, we can characterize the reduced Laurent norm unit ball as the set of points in the intersection half-spaces $H_{\alpha}$, one for each of the vertices of $\tilde{\Ne}(f))$.  
\begin{cor} With $f$ and $\Ne(f)$ as above, the reduced Laurent norm unit ball $ \tilde{\mathcal{B}}_f$ is an intersection taken over the vertices of $\Ne(f)$ of half-spaces\begin{displaymath} H_{\alpha}=\{ \phi \in (\mathbb{R}^m)^* \mid | \phi(\alpha)| \leq 1/2 \}\end{displaymath} so that \begin{equation}\nonumber  \tilde{\mathcal{B}}_f = \bigcap_{\alpha \in \ver(\tilde{\Ne}(f))}H_{\alpha}. \end{equation} \end{cor}
\section{The Alexander norm}
Before defining the Alexander norm and listing some of its properties as a special case of the Laurent norm, we give a brief historical summary of why this norm is important in the study of the topology of 3-manifolds, along with some recent results pertaining to this norm.
\subsection{Historical motivation for introducing the Alexander norm} The originator of the Alexander norm states in \cite{McM} that one motivation for introducing this norm in the late 1990s is the classical result from knot theory relating the degree of the Alexander polynomial of the knot to the genus $g$ of the knot; \begin{equation}\label{E:Tork} \deg \triangle(t)  \leq 2 g . \end{equation} This result establishes a relation between the degree of the Alexander polynomial and the Euler characteristic of the Seifert surface of the complement of the knot. This can be shown as follows. For any compact, connected, oriented surface $S$ of genus $g$ with $r$ holes, the Euler characteristic $\chi(S)$ for $S$ is given by the formula $\chi(S) = - 2g - r +2$. For the Seifert surface $S$ of the complement of a knot, $\chi(S) = -2g+1$. In terms of the Euler characteristic, the inequality of Equation \eqref{E:Tork} becomes \begin{equation}\nonumber \deg \triangle(t) - 1 \leq | \chi(S) |=2g-1 .\end{equation}Since the well-known Thurston norm is determined by the Euler characteristics of such surfaces, this suggests a deep connection between the Thurston norm, introduced in 1986 in \cite{Thu}, and the Alexander polynomial. We note that by Proposition \ref{P:Mman}, $\deg \triangle(t) = \| \phi\|_A$ for the cohomology class $\phi = (1)$.

In 1953, Torres \cite{Tor} proved that the multivariable Alexander polynomial $\triangle(t_1, \dots ,t_r)$ obtained from the Wirtinger presentation of an $r$-component link can be used to obtain an analogous relation for links to this result given by Equation \eqref{E:Tork} for knots. He defines a single variable Alexander polynomial by setting each variable $t_i$ equal to $t$ for $i=1, \dots ,r$ in the multivariable Alexander polynomial, and shows this single variable polynomial satisfies the inequality \begin{equation} \label{E:Tor} \deg \triangle(t, \dots ,t) \leq | \chi(S)| = 2g+r-2.\end{equation} We note that by Proposition \ref{P:Mman}, $\deg \triangle(t, \dots ,t) = \| \phi \|_A$ for the cohomology class $\phi=(1, \dots ,1)$.  In \cite{McM}, McMullen shows that the Alexander norm $\| \phi \|_A$ and Thurston norm $\| \phi \|_T$ satisfy a similar inequality for 3-manifolds with first Betti number $b_1 \geq 2$; $\| \phi \|_A \leq \| \phi \|_T$. In \cite{Tor}, Torres also showed that for alternating links the inequality of Equation \eqref{E:Tor} became an equality. In 2006, Ozv\'ath and Szaba\'o show in \cite{Osz} that the Alexander and Thurston norms coincide for alternating links using a proof based on Floer homology.

\begin{rem} The Alexander polynomial of a knot has an extra factor of $(t-1)$ in it as opposed to the Alexander polynomial of a link with two or more components.   For this reason, the degree of $\triangle(t)$ which appears in inequality of Equation \eqref{E:Tor} for a link has to be shifted by one in comparing it to the analogous result of Equation \eqref{E:Tork} for a knot. \end{rem} 

The Thurston norm has been found especially useful in the study of fibrations of a compact, oriented 3-manifold over the circle (see \cite{Thu}). However, it is often difficult to calculate since it depends on finding the Euler characteristics of the embedded, compact, connected, orientable surfaces in the 3-manifold. In contrast, the Alexander norm being determined directly from the Alexander polynomial of the 3-manifold, can always be determined by a straight-forward calculation. McMullen found that, in many cases, the Alexander norm coincides with the Thurston norm.  However, there are exceptional cases when the two norms do not coincide. An example of a two-component link for which this is the case was found by Dunfield \cite{Dun}.

We show in \cite{Lon1} that for a large class of links in homology 3-spheres that the Alexander and Thurston norm coincide. This class of links consists of the links which can be built up by splicing together link components which are Seifert links in homology 3-spheres. 
\subsection{The Alexander norm as a special case of the Laurent norm}Our definition of the Laurent norm is a generalization of the definition given by McMullen in \cite{McM} for the Alexander norm which is the Laurent norm in the case that the polynomial $f$ is also an Alexander polynomial $\triangle$ of a connected, compact, orientable 3-manifold whose boundary (if any) is a union of tori. McMullen uses the notation $\| \phi \|_A$ for the Alexander norm of cohomology class $\phi$ so that $\| \phi \|_A: = \| \phi \|_{\triangle}$. 
\begin{defn} An Alexander polynomial $\triangle$ determines a semi-norm called the \textit{Alexander norm} which is equivalent to the Laurent norm as defined in Definition \ref{D:fn} in the case that the Laurent polynomial $f$ is also an Alexander polynomial $\triangle$. The Alexander norm for $\triangle$ of $\phi$ is denoted as $\| \phi \|_A$ or  $\| \phi \|_{\triangle}$.  \end{defn}

In the case of knots and links, we can assume that the Alexander polynomial is determined by the complement $X$ of the knot or link which is just a finite CW complex. Then the Alexander norm is defined on not simply vectors $\phi \in (\mathbb{R}^n)^*$ but cohomology classes with real coefficients $\phi \in H^1(X;\mathbb{R})$. Further, the exponents of the Alexander polynomial which we have denoted $\alpha$ are homology classes $\alpha \in  H_1(X;\mathbb{Z})$. 

More generally, if the finite CW complex from which we obtain the Alexander polynomial corresponds to some $3$-manifold $M$, then we can define the cohomology classes on which the Alexander norm acts as follows. First we can assume that the first fundamental group $\pi(M)$ is finitely generated. The first homology group in this case is not just a free abelian group as for knots and links but it may have torsion. If we define the free abelian group $H$ with first Betti number $b_1$ as $H = H_1(M)/\Tor(H_1(M)) = \mathbb{Z}^{b_1}$, then the exponents of the Alexander polynomial lie in $H$. Even more, we have that $H^1(M; \mathbb{R}) = \Hom_{\mathbb{R}}(H_1(M),\mathbb{R}) = \Hom_{\mathbb{R}}(H,\mathbb{R}) = (H \otimes \mathbb{R})^* = (\mathbb{R}^{b_1})^*$ since contributions from the torsion part $\Tor(H_1(M))$ of $H_1(M)$ go to zero in passing to duality. Thus for 3-manifold $M$, the Alexander norm is a norm on vectors $\phi \in  H^1(M;\mathbb{R})= (\mathbb{R}^{b_1})^*$.

Since the Alexander norm for an Alexander polynomial $\triangle$ is just a special case of the Laurent norm for a polynomial which is also an Alexander polynomial, all our results above for the Laurent norm can be applied to directly to the Alexander norm.  Further, since Alexander polynomials are symmetric polynomials, we can also apply our results for the Laurent norms of symmetric polynomials to the Alexander norms.  Since an Alexander polynomial $\triangle$ is in the polynomial ring $ \mathbb{Z}[t_1^{\pm 1}, \dots ,t_{b_1}^{\pm 1}]$, the first Betti number $b_1$ of the 3-manifold replaces the integer $n$ in our discussion for the Laurent norm.  In addition, if we define $b_e$ as the \textit{essential first Betti number} which is equal to the dimension of the Newton polyhedron $\Ne(\triangle)$ of Alexander polynomial $\triangle$, then $b_e$ replaces the integer $m$ in our previous discussion of the Laurent norm.
\begin{rem} Essential variables of an Alexander polynomial have already been used in Dimca, Papadima and Suciu \cite{Dim}. They find that a Seifert link has only one essential variable, regardless of the number of link components, and conclude that the Newton polyhedron of a Seifert link must be a line segment.\end{rem}

\section{Two examples of Alexander norm calculations}\label{S:E}
In this section, we consider two examples of Alexander norm calculations for links. In the first, the number of essential variables is equal to the first Betti number of the link complement; $b_e = b_1 = 3$.  In the second, the number of essential variables is only two while the first Betti number is six; $b_e=2$ and $b_1=6$. 
\begin{ex}[The Borromean rings] The Borromean rings (see \cite{Burd2}) is a classic three-component link whose link complement has a hyperbolic geometry.  Its Alexander polynomial is 
\begin{eqnarray}\nonumber \triangle(t_1,t_2,t_3)  &= &(t_1 -1)(t_2 - 1)(t_3 - 1)\\ \nonumber
 &=& -1+t_1+t_2+t_3-t_1t_2-t_1t_3-t_2t_3+t_1t_2t_3.  \end{eqnarray} 
The vertex set of the Newton polyhedron, which coincides with the exponent set of the Alexander polynomial, is given by 
\begin{eqnarray}\nonumber   \ver(\Ne(\triangle))
 =  \lbrace (0,0,0),(1,0,0),(0,1,0),(0,0,1),\\ \nonumber (1,1,0),(1,0,1),(0,1,1),(1,1,1)\rbrace.    \end{eqnarray}
The Newton polyhedron, which is the convex hull of these points, is an octahedron.  Since $\Ne(\triangle)$ is three-dimensional, the space spanned by the essential variables is $\mathbb{R}^3$ and there are no inessential variables. Hence, by Theorem \ref{T:Rln}, we know that $\mathcal{\tilde{B}}_{\triangle} = \mathcal{B}_{\triangle}$. 
Writing the Alexander polynomial in terms of its irreducible factors, we have that $\triangle = f_1  \cdots  f_3$ where $f_i = t_i-1, i=1,\dots ,3$. Applying the Laurent norm decomposition formula, Equation \eqref{E:f1}, we obtain \begin{eqnarray}\label{E:Br1} \| \phi \|_{\triangle}&=& \| \phi \|_{(t_1-1)(t_2-1)(t_3-1)}
\\ \nonumber & =& \| \phi \|_{(t_1-1)} + \| \phi \|_{(t_2-1)} + \| \phi \|_{(t_3-1)} \\ \nonumber
& =& | \phi_1 | + | \phi_2 | + | \phi_3 |.  \end{eqnarray}
 By Lemma \ref{L:CL}, the vertices of the Alexander norm unit ball  can found by setting two out of the three absolute values in Equation \eqref{E:Br1} to zero and solving the equality $\| \phi \|_{(t_1-1)(t_2-1)(t_3-1)}=1$.  Since the Alexander norm unit ball is symmetric with respect to the origin, we also know that if $v$ is a vertex then $-v$ is also a vertex\begin{enumerate} \item $| \phi_1| = 0$, $| \phi_2| = 0$ and $|\phi_3| = 1 \Rightarrow (0,0,1)$ and $(0,0,-1)$ are vertices.
\item $| \phi_1| = 0$, $|\phi_3| =0$ and $|\phi_2| = 1\Rightarrow (0,1,0)$ and  $(0,-1,0)$ are vertices.
\item $| \phi_2| = 0$, $| \phi_3 |=0$, and $|\phi_1| = 1 \Rightarrow (1,0,0)$ and $(-1,0,0)$ are vertices.
\end{enumerate} The full set of vertices is 
\begin{equation}\nonumber \ver(\mathcal{B}_{\triangle}) = \lbrace (1,0,0),(-1,0,0),(0,1,0),(0,-1,0), (0,0,1),(0,0,-1) \rbrace. \end{equation}  
We thus find that the Alexander norm unit ball is the unit cube, which is the same result obtained in \cite{McM} by McMullen.\end{ex}
\begin{ex}[The great circle link $L(\mathcal{A})(321456)$]
We use as an example of an Alexander polynomial with an unbounded Alexander norm unit ball the great circle link $L(\mathcal{A})(321456)$. It is a 6-component link and is called a great circle link because each component traverses a circumference of $S^3$. This link is an example of an Eisenbud-Neumann graph link (see \cite{Eis}) since it can be formed by cabling on the unknot. Since it is a graph link in $S^3$, it can be constructed by splicing together Seifert links, each of which has a complement in $S^3$ that is Seifert-fibered. Links of this type are closely related to plane arrangements in $\mathbb{R}^4$, and a study of them can be found in Matei and Suciu \cite{Mat1}. The linking numbers between pairs of link components for this type of link are all $\pm 1$ and the self-linking number of each component is $0$.  This particular link is characterized by having all linking numbers $+1$ between pairs, except between any pair chosen from the first, second and third components.

The Alexander polynomial $\triangle$ for this link is given by \begin{equation} \nonumber \triangle = (t_1 \cdots t_6 - 1)^2(t_1^{-1} t_2^{-1}t_3^{-1} t_4 t_5t_6 - 1)^2.\end{equation} The convex hull of the exponents of this polynomial, which is the Newton polyhedron $\Ne(\triangle)$, is a two-dimensional polyhedron. This implies that there are two essential variables and four inessential variables. The Alexander norm can be found by letting $f_1 =(t_1 \dots t_6 - 1), f_2 = (t_1^{-1} t_2^{-1}t_3^{-1} t_4 t_5t_6 - 1), n_1 = 2,$ and  $n_2 = 2$ in our Laurent norm decomposition formula given by Equation \eqref{E:f1}. The result is \begin{eqnarray} \nonumber  \| \phi \|_{\triangle}&=&\| \phi\|_{(t_1 \cdots t_6 - 1)^2(t_1^{-1} t_2^{-1}t_3^{-1} t_4 t_5t_6 - 1)^2}\\ \nonumber &=&2 \| \phi\|_{(t_1 \cdots t_6 - 1)}+2 \| \phi\|_{(t_1^{-1} t_2^{-1}t_3^{-1} t_4 t_5t_6 - 1)}\\ \nonumber &=&2 | \sum_{i=1}^6 \phi_i| + 2| -\sum_{i=1}^3\phi_i + \sum_{i=4}^6 \phi_i|.  \end{eqnarray}  The set of points $\phi$ whose coordinates  satisfy $\sum_{i=1}^6 \phi_i=0$ and also $-\sum_{i=1}^3 \phi_i + \sum_{i=4}^6 \phi_i=0$ are non-zero vectors $\phi$ which have zero Alexander norm, $\| \phi \|_{\triangle}=0$. Hence, in this example, the Alexander norm is degenerate; it is only a semi-norm rather than a norm. By making the change of coordinates $\tilde{t}_1=t_1 + t_2 +t_3, \tilde{t}_2=t_4 +t_5 +t_6, \tilde{t}_i = t_i, (i=3,\dots , 6)$, we find that the Alexander polynomial in terms of essential coordinates is \begin{equation} \nonumber \triangle = (\tilde{t}_1 \tilde{t}_2 -1)^2(\tilde{t}_1^{-1} \tilde{t}_2 - 1)^2.\end{equation} In this coordinate system, it is easy to determine that the reduced Newton polyhedron, which is the convex hull of the exponents, is a two-dimensional diamond. Thus, there are two essential coordinates $\tilde{t}_1$ and $\tilde{t}_2$ and four inessential coordinates $\tilde{t}_3, \dots , \tilde{t}_6$. Passing to the dual coordinate system, there are two essential coordinates $\tilde{\phi}_1$ and $\tilde{\phi}_2$ which completely determine the Alexander norm. \begin{equation} \label{E:L1} \| \phi \|_{\triangle}=\| \phi \|_{(\tilde{t}_1 \tilde{t}_2 -1)^2(\tilde{t}_1^{-1} \tilde{t}_2 - 1)^2} = 2|\tilde{\phi}_1 + \tilde{\phi}_2| +2 | -\tilde{\phi}_1 + \tilde{\phi}_2|.\end{equation}   In essential coordinates, the Alexander norm vanishes only when $\phi=0$, so it is non-degenerate in the space of these coordinates in agreement with our Theorem \ref{T:Nd}.   Since there are four inessential coordinates, by Theorem \ref{T:Rln} the Alexander norm unit ball $\mathcal{B}_{\triangle}$ and the reduced Alexander norm unit ball $\tilde{\mathcal{B}}_{\triangle}$ are related by $  \mathcal{B}_{\triangle} = \tilde{\mathcal{B}}_{\triangle} \times (\mathbb{R}^4)^*$. To find the vertices of the reduced Alexander norm unit ball we apply Lemma \ref{L:CL} to Equation \eqref{E:L1}:  \begin{enumerate} \item $2|\tilde{\phi}_1 + \tilde{\phi}_2| = 0$  and  $2|\tilde{\phi}_1 - \tilde{\phi}_2| = 1$ implies that $(\pm 1/4, \mp1/4)$ are vertices.  
\item $2|\tilde{\phi}_1 - \tilde{\phi}_2| = 0$  and  $2|\tilde{\phi}_1 + \tilde{\phi}_2| = 1$ implies that $(\pm 1/4, \pm 1/4)$ are vertices.\end{enumerate} The full set of vertices is \begin{equation}\nonumber \ver(\mathcal{\tilde{B}}_{\triangle}) = \lbrace (1/4, 1/4),(-1/4,-1/4),(1/4,-1/4),(-1/4,1/4) \rbrace. \end{equation}  The reduced Alexander norm unit ball $\tilde{\mathcal{B}}_{\triangle}$, which is the convex hull of these vertices, is a square.
\end{ex}

\textbf{Acknowledgements }

This article is based on Ph.D. thesis research conducted at Northeastern University under the direction of Alexander~I.~Suciu. Further details can be found in Long \cite{Lont}, \cite{Lon1}.

\end{document}